\documentclass{article}
\usepackage{cmap,amsmath,amsthm,amssymb}
\usepackage[utf8]{inputenc}
\usepackage[english]{babel}
\usepackage{pstricks}
\usepackage[dvips]{graphicx}
\usepackage{pst-3dplot}
\usepackage[all,cmtip]{xy}
\usepackage{tikz}
\usepackage{pst-all}
\usepackage{pstricks-add}
\usepackage{enumitem}
\usepackage{pst-solides3d}
\usepackage{subfigure}
\usepackage{color, graphics}
\usepackage{float}

\newtheorem{thm}{Theorem}

\newtheorem*{thm*}{Theorem}

\title{On differential equations of integrable billiard tables}
\author{Vladimir Dragovi\'c$^1$ and Andrey E. Mironov$^2$}
\date{}

\begin{document}
\sloppy
\maketitle

\footnotetext[1]{Department of Mathematical Sciences, The University
	of Texas at Dallas, Richardson TX,
	USA; Mathematical Institute SANU,
	Belgrade, Serbia.  E-mail: {\tt
		Vladimir.Dragovic@utdallas.edu}}

\footnotetext[2]{Novosibirsk State University,
Novosibirsk,
Russia;
Sobolev Institute
of Mathematics of the Siberian Branch of the Russian Academy of Sciences,
Novosibirsk,
Russia. E-mail: {\tt mironov@math.nsc.ru}}

\begin{abstract}
We introduce a method to find differential equations for functions which define tables, such that associated billiard systems admit a local first integral. We illustrate this method in three situations: the case of (locally) integrable wire billiards, for finding surfaces in ${\mathbb R}^3$ with a first integral of degree one in velocities, and for finding a piece-wise smooth surface in ${\mathbb R}^3$ homeomorphic to a torus, being a table of an integrable billiard.
\end{abstract}

\

{\it Dedicated to Misha Bialy, a colleague and friend, on the occasion of his 60-th anniversary.}

\

\

\section{Introduction}

Mathematical billiards provide a very important class of dynamical systems, see for example \cite{KT, Tab}. The integrability of such systems
has been intensively studied from various perspectives, see e.g. \cite{bol, KT, DR, DR1, ADSK, KS, bm, bm1, glu1, ShTre} and references therein.

 In this paper we obtain differential equations on tables for different types of billiards which  admit  first integrals polynomial in components of the velocity vector. Further, for brevity,  we will refer to such first  integrals just as polynomial integrals.  We apply our method first  to wire billiards introduced in \cite{bmt}. We show that in the case of integrable wire billiards constructed in \cite{bmt}, there is a polynomial integral of degree one (Theorem~1, see bellow). Then, we find surfaces in ${\mathbb R}^3$ which admit a first integral of degree one  (Theorem 2). We also find piece-wise smooth surfaces in ${\mathbb R}^3$ homeomorphic to a torus, each being a table for an integrable billiard.

 To begin with, we explain our method in the case of Birkhoff billiards \cite{KT, bol} in the plane.
 Let $\gamma(t)=(\gamma_1(t),\gamma_2(t))$ be a convex smooth curve parametrized by the arc length, $|\dot{\gamma}|=1$. A particle moves along straight lines inside the table, bounded by the curve $\gamma$. We assume that the speed of the particle is equal to unity there, i.e. $|v|=1$, where $v=(v_1,v_2)$ is the velocity vector of the particle. When  the particle reaches the boundary of the table, defined by a curve $\gamma$, it is reflected according to the geometric optics law. Any function of the form
 $F(M,v_1,v_2),$ where $M=xv_2-yv_1$ is constant along trajectories between reflections. On the other hand, any function $G(S,t)$, where  $S=(v,\dot{\gamma})=\dot{\gamma}_1v_1+\dot{\gamma}_2v_2$, is constant at the moment of a reflection. Hence, if the identity
 $$
  F(\gamma_1(t)v_2-\gamma_2(t)v_1,v_1,v_2)=G((v,\dot{\gamma}),t)\eqno{(1)}
 $$
 holds for all $v,|v|=1$, then the Birkhoff billiard is integrable within the domain bounded by $\gamma(t)$. Let us consider some examples.

 \noindent{\bf Example 1 (circle).} Let us consider the case when $F$ is a polynomial integral of degree one
 $$
   F=xv_2-yv_1+av_1+bv_2,\ a,b\in{\mathbb R}.
 $$
 In this case the condition (1) has the form
 $$
  \gamma_1(t)v_2-\gamma_2(t)v_1+av_1+bv_2=h(t)(\dot{\gamma}_1(t)v_1+\dot{\gamma}_2(t)v_2),
 $$
 From here we obtain the system of differential equations
 $$
 \begin{aligned}
  a-\gamma_2(t)-h(t)\dot{\gamma}_1(t)&=0,\\
 b+\gamma_1(t)-h(t)\dot{\gamma}_2(t)&=0.
 \end{aligned}
$$
 Hence
 $$
  a\dot{\gamma}_2(t)-\dot{\gamma}_2(t)\gamma_2(t)-b\dot{\gamma}_1(t)-\dot{\gamma}_1(t)\gamma_1(t)=0,
$$
or
$$
 (\gamma_1(t)+b)^2+(\gamma_2(t)-a)^2=R^2.
$$
We obtain that the curve $\gamma$ is a circle.

 \noindent{\bf Example 2 (parabola).} We consider the case when $F$ has the form
 $$
 F=M v_2+\lambda v_1^2,\ \lambda\in{\mathbb R}.
 $$
 The equation (1) has the form
 $$
 (\gamma_1(t)v_2-\gamma_2(t)v_1)v_2+\lambda v_1^2=h_1(t)(\dot{\gamma}_1(t)v_1+\dot{\gamma}_2(t)v_2)^2+h_2(t).
 $$
 Using $v_2^2=1-v_1^2$, we get the system of equations
 $$
 \begin{aligned}
 \lambda-\gamma_1(t)-h_1(t)(\dot{\gamma}_1^2(t)-\dot{\gamma}_2^2(t))&=0,\\
 \gamma_2(t)+2h_1(t)\dot{\gamma}_1(t)\dot{\gamma}_2(t)&=0,\\
 h_2(t)-\gamma_1(t)+h_1(t)\dot{\gamma}_2^2(t)&=0.
 \end{aligned}
 $$
 Hence
 $$
 \lambda-\gamma_1(t)+\frac{\gamma_2(t)}{2\dot{\gamma}_1(t)\dot{\gamma}_2(t)}(\dot{\gamma}_1^2(t)-\dot{\gamma}_2^2(t))=0.
 $$
 One can check by a direct calculation that the solution of this equation is a conic, a parabola:
 $$
 p^2-2p\lambda+2p\gamma_1(t)=\gamma_2^2(t), \ p\in{\mathbb R}.
 $$

 \noindent{\bf Example 3 (ellipse, hyperbola).} Let us consider the case when $F$ has the form
 $$
  F=M^2+\lambda v_1^2,\ \lambda\in{\mathbb R}.
$$
The equation (1) gives us
$$
 (\gamma_1(t)v_2-\gamma_2(t)v_1)^2+\lambda v_1^2=h_1(t)(\dot{\gamma}_1(t)v_1+\dot{\gamma}_2(t)v_2)^2+h_2(t).
 $$
 Using $v_2^2=1-v_1^2$ we get the system of equations
 $$
 \begin{aligned}
  \lambda-\gamma_1^2(t)+\gamma_2^2(t)-h_1(t)(\dot{\gamma}_1^2(t)-\dot{\gamma}_2^2(t))&=0,\\
 \gamma_1(t)\gamma_2(t)+h_1(t)\dot{\gamma}_1(t)\dot{\gamma}_2(t)&=0,\\
 h_2(t)-\gamma_1^2(t)+h_1(t)\dot{\gamma}_2^2(t)&=0.
 \end{aligned}
$$
Thus,
$$
 \lambda-\gamma_1^2(t)+\gamma_2^2(t)+\frac{\gamma_1(t)\gamma_2(t)}{\dot{\gamma}_1(t)\dot{\gamma}_2(t)}(\dot{\gamma}_1^2(t)-\dot{\gamma}_2^2(t))=0.
$$
One can easily see that the solution of this equation is a conic, a hyperbola or an elipse:
$$
 \frac{\gamma_1^2(t)}{a^2}+\frac{\gamma_2^2(t)}{a^2-\lambda}=1.
$$

  \section{Wire billiards}

 Let $\gamma$ be a smooth curve  in ${\mathbb R}^n$ which will play a role of a wire defining a wire billiard, as follows.  The chord
 $\gamma(s)\gamma(t)$ reflects to a chord $\gamma(t)\gamma(s_1)$ if the angles between
 the chords and the tangent vector
 to $\gamma$ at $\gamma(t)$ are equal.
 Generally, a point $\gamma(s_1)$ is  not uniquely defined.
 In \cite{bmt} an example of integrable wire billiards is given.

 Let $A\in so(n), \gamma_0\in{\mathbb R}^n$. Consider the curve
$$
 	\gamma(t)=e^{At}\gamma_0\subset{\mathbb R}^n.
 $$
 It turns out that the angles between the segments and $\gamma$ are preserved under the reflection. Thus, this is an example of a completely integrable wire billiard in ${\mathbb R}^n$.

 \noindent{\bf Example 4.} The wire billiard defined by the curve in ${\mathbb R}^4={\mathbb C}^2$ given by the formula
 $$
 \gamma(t)=(ae^{ikt},be^{imt}), a,b>0, t\in[0,2\pi],k,m\in{\mathbb N}
 $$
 is integrable. The curve is a toric knot in $S^3\subset{\mathbb R}^4.$ Although the integrability of this system was shown in  \cite{bmt}, the form of the first integral was not studied there.

 In this paper we show that the wire billiard defined by the curve $\gamma(t)$ admits a first integral, which is a polynomial of degree one.

 Let $M_{ij}=v_jx_i-v_ix_j, 1\leq i<j\leq n,$ where $v=(v_1,\dots,v_n)$ is the velocity vector of the particle, $|v|=1$, and $x=(x_1,\dots,x_n)$  are the coordinates of the particle.

 \begin{thm} The wire billiard defined by the curve $\gamma=e^{At}\gamma_0\subset{\mathbb R}^n$ admits the first integral
 	$$
 	F=\sum_{i<j}a_{ij}(v_jx_i-v_ix_j),
 	$$
 	where $a_{ij}$ are the components of the matrix $A$.
 \end{thm}

\noindent {\bf Proof.  }
 A function of the form
 $$
 F(M_{ij},v)=F(M_{12},\dots,M_{n-1n},v_1,\dots,v_n)
 $$
 is constant along the motion of the particle between reflections.
 Let us consider a point $\gamma(t)$, $|\dot{\gamma} (t)|=1.$
 Any function of the form
 $$
  G(\cos\varphi,t)=G((v,\dot{\gamma}),t)=G(v_1\dot{\gamma}_1+\dots+v_n\dot{\gamma}_n,t),
 $$
 where $\varphi$ is the angle between the velocity vector $v$ and $\dot{\gamma}$,
 is invariant under the reflection. Thus, if the identity
 $$
 F(v_j\gamma_i(t)-v_i\gamma_j(t),v)=G(v_1\dot{\gamma}_1+\dots+v_n\dot{\gamma}_n,t)\eqno{(2)}
 $$
 is satisfied for all $v$, then the wire billiard is integrable.

 Let us consider the simplest case when $F$, and thus $G$, is a polynomial integral of degree one. We have
 $$
 \sum_{i<j}a_{ij}(v_j\gamma_i(t)-v_i\gamma_j(t))+b_1v_1+\dots b_nv_n=\dot{\gamma_1}v_1+\dots+\dot{\gamma_n}  v_n.
 $$
 This condition is equivalent to the system of equations
 $$
 \dot{\gamma}=A\gamma+b,\ b=(b_1,\dots,b_n)^{\top},\ A=(a_{ij}),
 $$
 where $A+A^{\top}=0$. In the case $b=0$ we have
 $$
 \dot{\gamma}=A\gamma,\quad
 \gamma(t)=e^{At}\gamma_0.
 $$
 \noindent{\bf Example 5.} It would be very interesting to construct an integrable wire billiard for a closed curve in
 ${\mathbb R}^3$. In Example 4 the wire is a non-closed curve. We consider this case in more detail. After an appropriate rotation and a shift of the coordinate system one can assume that the integral has the form
 $$
  F=yv_1-xv_2+av_3,
 $$
 then (2) gets the form
 $$
  \gamma_2(t)v_1-\gamma_1(t)v_2+av_3=h(t)(v_1\dot{\gamma}_1+v_2\dot{\gamma}_2+v_3\dot{\gamma}_3).
 $$
 From here, we get the system
 $$
 \begin{aligned}
  \gamma_2(t)-h(t)\dot{\gamma}_1&=0,\\
  \gamma_1(t)+h(t)\dot{\gamma}_2&=0,\\
  a-h(t)\dot{\gamma}_3&=0.
  \end{aligned}
 $$
 This system has a solution
 $$
  \gamma_1(t)=R\sin(t),\  \gamma_2(t)=R\cos(t),\ \gamma_3(t)=at,\ h(t)=1.
 $$
 The solution presents a circular spiral.

 \section{Surfaces with a first billiard integrals of degree one}
 In this section we consider surfaces in ${\mathbb R}^3$ which admit a first (local) billiard integral of degree one in components of the velocity vector.

 Let $\Sigma\subset{\mathbb R}^3$ be a surface given as an image of
 $$
 r:U\rightarrow{\mathbb R}^3,\quad r(u)=(r^1(u),r^2(u),r^3(u)), u\in U.
 $$
 We introduce the functions
 $$
 \begin{aligned}
  M_1&=r^2(u)v_3-r^3(u)v_2,\\
  M_2&=r^3(u)v_1-r^1(u)v_3,\\
  M_3&=r^1(u)v_2-r^2(u)v_1.
  \end{aligned}
 $$
 A possible first billiard integral of degree one has the form
 $$
 F=a_1 M_1+a_2 M_2+a_3 M_3+b_1v_1+b_2v_2+b_3v_3.
 $$
 The function $F$ is constant along the motion between reflections. Now we find a condition for $F$ to be also preserved under the reflections.

 Making an appropriate rotation of the orhogonal coordinate system, we get $F$ in the following form
 $$
 F= \alpha M_3+b_1v_1+b_2v_2+b_3v_3.
 $$

 By applying an appropriate shift of the coordinate system, without loss of generality we get that
 $$
 F= \alpha M_3+\beta v_3.
 $$
 Denote by
 $$
 \begin{aligned}
 S_1&=(v,r_{u_1})=v_1r^1_{u_1}+v_2r^2_{u_1}+v_3r^3_{u_1},\\
 S_2&=(v,r_{u_2})=v_1r^1_{u_2}+v_2r^2_{u_2}+v_3r^3_{u_2}.
 \end{aligned}
 $$
 Any function of the form
 $$
 G(S_1,S_2,u_1,u_2)
 $$
 is constant under the reflection at the point of the surface with coordinates $(u_1, u_2)$.
 Thus, we have an identity
 $$
 \alpha M_3+\beta v_3=h_1(u) S_1+h_2(u)S_2,
 $$
 where $h_1(u), h_2(u)$ are some functions.
 From the last identity we have
 $$
 \begin{aligned}
 \alpha r^2+h_2r^1_{u_2}+h_1r^1_{u_1}&=0,\\
 \alpha r^1-h_2r^2_{u_2}-h_1r^2_{u_1}&=0,\\
 \beta-h_2r^3_{u_2}-h_1r^3_{u_1}&=0.
 \end{aligned}
 $$
 Hence
 $$
 h_1(u)=\frac{\alpha(r^1r^1_{u_2}+r^2r^2_{u_2})}{r^1_{u_2}r^2_{u_1}-r^1_{u_1}r^2_{u_2}},\quad
 h_2(u)=\frac{\alpha (r^1r^1_{u_1}+r^2r^2_{u_1})}{r^1_{u_1}r^2_{u_2}-r^1_{u_2}r^2_{u_1}}.
 $$
 We put
 $$
 r^1=u_1\qquad r^2=u_2.
 $$
 Then we get the equation
 $$
 \beta-\alpha u_1r^3_{u_2}+\alpha u_2r^3_{u_1}=0.
 $$
 This equation has a solution
 $$
  r^3=-\frac{\beta}{\alpha}ArcTan(\frac{u_1}{u_2})+f(u_1^2+u_2^2),
 $$
 where $f(x)$ is a function of one variable.
 We obtain
 \begin{thm}  A surface in ${\mathbb R}^3$ parametrized by
 	$$
 	 r(u)=\left (u_1,u_2,-\frac{\beta}{\alpha}ArcTan(\frac{u_1}{u_2})+f(u_1^2+u_2^2)\right)
 	$$
 	admits a first (local) billiard integral of the form
 	$$
 	F=\alpha (u_1v_2-u_2v_1)+\beta v_3.
 	$$
 \end{thm}

 \section{Integrable biliards inside piecewise smooth surfeces homeomorphic to a torus}

 In this section we construct a piece-wise smooth surface homeomorphic to a torus with two independent first  billiard integrals.
 We assume that
 $$
 r^1(u)=u_1,\ r^2(u)=u_2,\ r^3(u)=f(u_1^2+u_2^2).
 $$
 Then, by applying  Theorem 2 to the case $\beta=0$, we see that this surface admits a first billiard integral of degree one
 $$
 F_1=M_3.
 $$
 We assume that there is an additional first integral of degree two of the form
 $$
 F_2=a(M_1^2+M_2^2)+b(M_1v_2-M_2v_1)+c(v_1^2+v_2^2),\ a,b,c\in {\mathbb R}.
 $$
 The condition that $F_2$ is preserved under the reflection implies that $F_2$ has the form
 $$
 F_2=h_{11}(u)S_1^2+h_{12}(u)S_1S_2+h_{22}(u)S_{2}^2+h(u),
 $$
 where $h_{11}(u), h_{12}(u), h_{22}(u), h(u)$ are some functions of one variable.
 From this condition and recalling that $v_1^2+v_2^2+v_3^2=1$, we obtain
 $$
 h_{11}(u)=\frac{b-2af(u_1^2+u_2^2)}{4f'(u_1^2+u_2^2)},\
 h_{12}(u)=0,\
 h_{22}(u)=\frac{b-2af(u_1^2+u_2^2)}{4f'(u_1^2+u_2^2)},
 $$
 and
 $$
 h=a(u_1^2+u_2^2)-4h_{11}(u)(u_1^2+u_2^2)f'(u_1^2+u_2^2),
 $$
 together with
 $$
 h_{11}(u)(1-4f'^2(u)(u_1^2+u_2^2))=(a-b)f(u)+(c-a(u_1^2+u_2^2))
 $$
 Thus $f$ satisfies an ordinary differential equation
 \begin{equation}\label{eq:ode}
 \mathcal F(t, f, f'):=b+4(at-c) f'-4af^2f'-4btf'^2-f(2a-4bf'-8atf'^2)=0,
 \end{equation}
 where we substitute $t=u_1^2+u_2^2$.
 We will discuss two types of these equations \eqref{eq:ode} depending on the conditions on the parameter $a$ and study the corresponding solutions

 \begin{figure}[h]
 \begin{minipage}{0.49\textwidth}
\includegraphics[width=6cm,height=6cm]{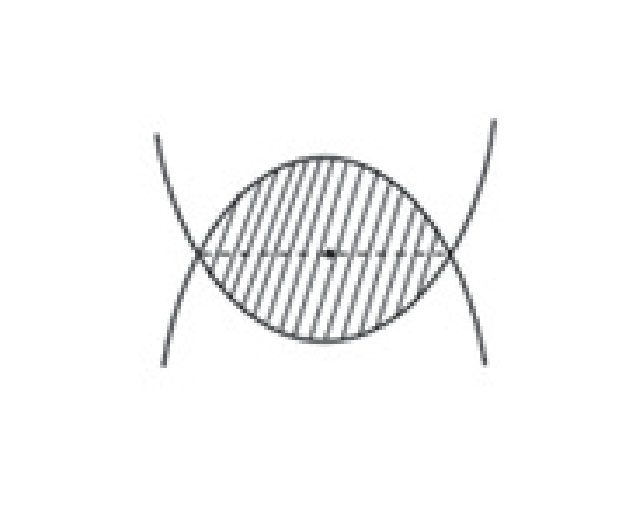}
\caption{A region bounded by two confocal parabolas.}
\label{fig:parabolas}
\end{minipage}
\begin{minipage}{0.49\textwidth}
\includegraphics[width=6cm,height=6cm]{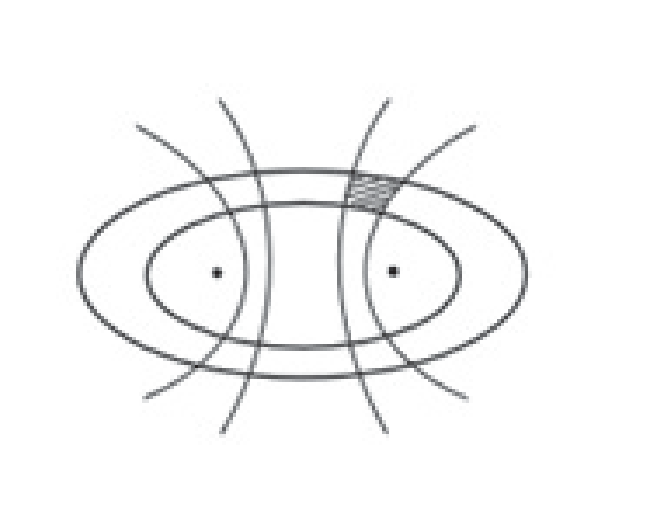}
\caption{A tetragon made of four confocal conics, two ellipses and two hyperbolas.}
\label{fig:ellipses}
\end{minipage}
\end{figure}

 If $a=0$, then, the function
 $$
 f(t)=st+\frac{4sc-b}{4sb},\ s\in{\mathbb R}.
 $$
 depending on $s$ as a parameter, generates solutions of \eqref{eq:ode}.

 If $a\ne 0$, then $f(t)$ defined from the equation
 $$
 \left(f(t)-\frac{b}{2a}\right)^2+\mathcal At-s=0,
 $$
 where
 $$
 \mathcal A=\frac{4a^2s}{4a^2s-b^2+4ac},
 $$
 and $s$ is a parameter, generates solutions of \eqref{eq:ode}.

 Using the first type of solutions, we construct an integrable billiard for a closed  piecewise smooth surface homeomorphic to a sphere, see Fig. \ref{fig:parabolas}. We get a family of confocal parabolas
 $$
 z=sR^2+\frac{4sc-b}{4sb}.
 $$

 For example, take $b=2, c=1$. Then for $s=1$ and $s=-1/2$ we obtain two surfaces in ${\mathbb R}^3$
 $$
 z=x^2+y^2+\frac{1}{4},\ z=-\frac{x^2+y^2}{2}+1.
 $$
 These two surfaces form a piecewise smooth surface homeomorphic to a sphere. This piecewise surface defines a Birkhoff billiard which admits two first integrals.

 In the second case,  we have a one-parametric family of confocal conics:
  $$
  \frac{z^2}{s}+\frac{R^2}{s +\frac{4ac-b}{4a^2}}=1.
  $$

  We choose the values of the parameter $s$ to obtain two ellipses and two hyperbolas. This way, we obtain a ``tetragon'', two sides of which are pieces of the ellipses and two sides are pieces of the hyperbolas, see Fig. \ref{fig:ellipses}. After rotating this tetragon around the axis $z$, we obtain a piecewise smooth surface, homeomorphic to a torus, with two billiard integrals, one of degree one and another of degree two.

\subsection*{Acknowledgements}
We are grateful to Misha Bialy for very helpful discussions. We dedicate this work to his 60th anniversary and wish him many happy returns.
The research of AM has been partially supported by Russian Science Foundation (Grant No. 21--41--00018) and of VD by the Science Fund of Serbia (Grant Integrability and Extremal Problems in Mechanics, Geometry and
Combinatorics, MEGIC, Grant No. 7744592), the Ministry for Education, Science, and Technological Development of Serbia, and the Simons Foundation (Grant No. 854861).

\end{document}